\input amstex
\documentstyle{amsppt}
\pagewidth{6.5in}
\pageheight{9in}
\topmatter
\title
A Characterization of Reflexivity \\for Dual Banach Spaces
\endtitle
\author Javier H. Guachalla \endauthor
\affil Universidad Mayor de San Andr\'es\footnote{e-mail:jguachal\@entelnet.bo} \\La Paz - Bolivia \endaffil
\address La Paz Bolivia \endaddress
\email jguachal\@entelnet.bo \endemail
\date August 2005 \enddate
\subjclass FA - Functional Analisis\endsubjclass
\abstract
Given a Banach space. We show that its three times dual space can be written as a direct sum. Being one of the sumands null is a necessary and sufficient condition for the reflexivity of the dual space. We end with an application to the known problem of the relationship between the reflexivity of a Banach space and that of its dual space.
\endabstract
\endtopmatter
\document

\head Introduction \endhead

Let $X$ be a Banach Space. We denote $X^*,\ X^{**}, \dots $, the (topological) dual, bidual, $\dots$ spaces of $X$, and 
$$\align
\pi&:X\rightarrow X^{**} \\
\pi'&:X^*\rightarrow X^{***}
\endalign
$$ 
the canonical maps into their respective bidual spaces.

For $X$, $Y$ Banach spaces, and $T:X\rightarrow Y$ a bounded linear operator given, we denote $Ran(T)$ the range of $T$ and $Ker(T)$ the null space of $T$. Note that since $T$ is bounded, $Ker(T)$ is a closed subspace of $X$.

The adjoint operator of the bounded linear operator $T$ will be denoted $T^*$, which is defined by $T^*(\mu)=\mu\circ T$ for $\mu \in Y^*$ and is a bounded linear operator; thus for the map $\pi$, we have 
$$
\pi ^{*}(\eta )=\eta \circ \pi \qquad (for\ all\ \eta \in X^{***}). 
$$

\subhead {1. Definition}\endsubhead Let $X$ and $Y$ be Banach spaces, $T:X\rightarrow Y$ be a bounded linear operator. We will say that T is an isomorphism if it has a bounded inverse linear operator.

\subhead {2. Definition}\endsubhead Let $X$ be a Banach space and $M,$ $N$ closed supspaces of $X$. We will say that $X$ is the direct sum of $M$ and $N$, denoted $X=M\oplus N$, if $X=M+N$ and $M\cap N=\{0\}$ where $M+N=\{x+y:x\in M \text{ and }y\in N\}$.

We begin with a theorem which extends a result from algebra.

\proclaim{3. Theorem} Let $X,$ $Y,$ and $Z$ be Banach spaces, $T:X\rightarrow Y$ and $S:Y\rightarrow Z$ bounded linear operators. Let the composition $ST$ be an isomorphism. Then the subspace $Ran(T)$ is closed and 
$$Y=Ran(T)\oplus Ker(S).$$
\endproclaim

\demo{Proof}
Let us see that $Ran(T)$ is closed. Let $y \in \overline{Ran(T)}$. There exist a sequence $\{x_n\}$ in $X$ such that 
$$ y=\lim_{n \to \infty}T(x_n)$$
Since $ST$ is an isomorphism, there exists a unique $x \in X$ such that $ST(x)=S(y)$. Let us prove that $x_n \to x$. We have that $(ST)^{-1}$ is bounded,
$$
\gather
ST(x_n)\to ST(x)\qquad \text{and}\\
\Vert x_n - x\Vert \le \Vert (ST)^{-1}\Vert\,\Vert (ST)(x_n)-(ST)(x)\Vert
\endgather
$$ 
Hence $$x_n \to x$$
Therefore $\lim_{n\to \infty}T(x_n)=T(x)$. Which means $y=T(x)$, i.e. $Ran(T)$ is a closed subspace.

Since $Ker(S)$ is a closed subspace, let us verify the direct sum conditions. 
Let $y\in Y.$ Then $\exists\,!\,x\in X$ such that $(ST)(x)=S(y).$
Let us call $y^{\prime }=y-T(x)$. Then $y^{\prime}\in Ker(S).$ Therefore, we can write $y=T(x)+y^{\prime }$. Where $T(x)\in Ran(T)$ and $y^{\prime }\in Ker(S)$. Thus $Y=Ran(T)+Ker(S)$.

Let us see that $Ran(T)\cap Ker(S)=\{0\}$. Suppose $y\in Ran(T)\cap Ker(S)$. Therefore $y=T(x)$, for some $x\in X$ and $S(y)=0$. Thus $S(T(x))=S(y)=0.$
Therefore $x=0$, since $ST$ is an isomorphism. So $y=T(0)=0.$
Hence, according to the definition
$$
Y=Ran(T)\oplus Ker(S)
$$
\enddemo

\head The dual space and its canonical map\endhead

In this section, we apply the above result to the canonical map of the dual space of a Banach space $X$ and to the adjoint of the canonical map of $X$ itself.

\proclaim{4. Theorem.} Let $X$ and $X^*$ be a Banach space and its dual space. Consider the maps $\pi ^{\prime }:X^{*}\rightarrow X^{***}$ and the adjoint operator of the canonical map of $X$, $\pi ^{*}:X^{***}\rightarrow X^{*}$. Then $\pi ^{*}\circ \pi^{\prime }$ is the identity map on $X^*$, is an isomorphism, and 
$$
X^{***}=Ran(\pi^{\prime })\oplus Ker(\pi ^{*}).
$$
\endproclaim

\demo{Proof} Let $\lambda \in X^{*}.$ Then, by definition we have
$$
(\pi ^{*}\circ \pi ^{\prime})(\lambda )=\pi ^{*}(\pi ^{\prime }(\lambda ))=\pi ^{\prime }(\lambda )\circ \pi
$$

Evaluating in an arbitrary $x\in X.$ 
$$
\align
[(\pi ^{*}\circ \pi ^{\prime })(\lambda )](x)&=[\pi^{\prime }(\lambda )\circ \pi ](x)\\
&=\pi ^{\prime }(\lambda )(\pi
(x))\\
&=(\pi (x))(\lambda )=\lambda (x)\qquad \text{(for all $x\in X$)}
\endalign
$$
therefore 
$$
\pi ^{*}\circ \pi ^{\prime }(\lambda )=\lambda \qquad (\text{for all $\lambda \in X^{*}$)} 
$$
Then we conclude that $\pi ^{*}\circ \pi ^{\prime }=1_{X^{*}}$ the identity map in $X^{*}$. Since $X^{*}$ has the same topology in both maps, the identity map, $1_{X^{*}}$ is an isomorphism. Hence, by theorem 3
$$X^{***}=Ran(\pi^{\prime})\oplus Ker(\pi^*)$$
\enddemo

\proclaim{5. Corollary} Let $X$ be a Banach space. The dual space $X^{*}$ of $X$ is reflexive if and only if $Ker(\pi ^{*})=0$. Where $\pi ^{*}$ is the adjoint map of the canonical map $\pi $ of $X$.
\endproclaim

\head Application\endhead

Let $X$ be a Banach space. We prove $X$ is reflexive if and only if $X^{*}$ is reflexive.\par
To do so, let us consider the following notations. Given a subset $M$ of a Banach space. We denote $M^o$ the annhilitor of $M$, that is, the closed subspace of functionals of the respective dual space, whose restriction to $M$ are zero.\par
Similarly, we denote ${}^oN$, for a subset $N$ of a dual space; the closed subspace of elements of the space itself, that anhilite $N$, i.e. if $N\subset X^*$ then ${}^oN$ is the set $\{x\in X: \forall \,f\in N\ f(x)=0\}$, which in fact is a closed subspace.

Suppose first, that $X^{*}$ is reflexive. By Corollary 5 $Ker(\pi ^{*})=0.$
Then by \cite{T.4.6}, and since $Ran(\pi)$ is a closed subespace, being $\pi$ an isometry, we obtain
$$
\align
Ran(\pi)&={}^oKer(\pi ^{*})\\
&={}^o\{0\}\\
&=X^{**}
\endalign
$$
Which says that $X$ is reflexive.

Let us now suppose that $X$ is reflexive, by the later reference and $Ran(\pi)$ closed, we have
$$
\align
Ker(\pi ^{*})&=Ran(\pi)^o\\
&=(X^{**})^o\\
&=\{0\}
\endalign
$$
by an application of the Hahn-Banach theorem,  and therefore, Corollary 5 implies $X^{*}$ is reflexive.

\end{document}